\documentclass{amsart}
\usepackage{amscd,amssymb,amsxtra}
\usepackage[mathscr]{eucal}
\usepackage[all]{xy}

\usepackage[usenames]{color}
\usepackage{graphicx}

\CompileMatrices

\title{Equivariant multiplicative closure}

\author{M.~A.~Hill}
\address{Department of Mathematics \\ University of Virginia
\\Charlottesville, VA 22904}
\email{michael.a.hill@math.uva.edu}

\author{M.~J.~Hopkins}
\address{Department of Mathematics \\ Harvard University
\\Cambridge, MA 02138}
\email{mjh@math.harvard.edu}

\newtheorem{cor}[equation]{Corollary}
\newtheorem{lem}[equation]{Lemma}
\newtheorem{prop}[equation]{Proposition}
\newtheorem*{thm*}{Theorem}
\newtheorem*{cor*}{Corollary}
\newtheorem*{lem*}{Lemma}
\newtheorem*{prop*}{Proposition}

\ifx\undefined\theoremstyle
	\newtheorem{definition}[equation]{Definition}
	\newenvironment{defin}{\begin{definition}\rm}{\end{definition}}
	\newtheorem{conj}[equation]{Conjecture}
	\newtheorem{rem}[equation]{Remark}
	\newtheorem*{rem*}{Remark}
	\newtheorem{rems}[equation]{Remarks}

\else
	\theoremstyle{definition}

	\newtheorem{defin}[equation]{Definition}

	\theoremstyle{remark}

	\newtheorem*{rem*}{Remark}

\ifx\undefined\prem	\newtheorem{prem}{\bf Working Remark} \fi

\fi

\newtheorem{eg}[equation]{Example}
\newtheorem{ques}[equation]{Question}

\ifx\undefined\pf
\newenvironment{pf}{\bigskip{\em Proof:\/}}{\qed\medskip}
\newenvironment{pf*}[1]{\bigskip{\em #1:\/}}{\qed\medskip}
\fi

\ifx\undefined\numberwithin
\def\numberwithin#1#2{\makeatletter\@ifundefined{c@#1}{\@nocnterrr}{%
  \@ifundefined{c@#2}{\@nocnterr}{%
  \@addtoreset{#1}{#2}%
  \toks@\expandafter\expandafter\expandafter{\csname the#1\endcsname}%
  \expandafter\xdef\csname the#1\endcsname
    {\expandafter\noexpand\csname the#2\endcsname
     .\the\toks@}}}\makeatother}\fi

\ifx\undefined\qed\newcommand{\qed}{\hfil\rule{4pt}{6pt}\bigskip}\fi
\ifx\undefined\\operatorname\newcommand{\operatorname}[1]{\mathop{\mbox{\rm #1}}}\fi

\newcommand{\Z}{{\mathbb Z}}

\newcommand{\zerowidth}[1]{\hbox to 0pt{\hss$\displaystyle #1$\hss}}

\ifx\undefined\eqref\newcommand{\eqref}[1]{\rm (\ref{#1})}\fi

\newcounter{thmItem}

\newcommand{\thmItemref}[1]
 	 {{\rm \ref{#1})}}

\newenvironment{thmList}{\begin{list}%
{\rm \roman{thmItem})}{\usecounter{thmItem}
\setlength{\labelwidth}{2em}
\setlength{\itemindent}{2em}
\setlength{\leftmargin}{0pt}
\setlength{\listparindent}{0pt}
\setlength{\parsep}{0pt}
\setlength{\partopsep}{0pt}
\setlength{\itemsep}{\medskipamount}
\setlength{\topsep}{\medskipamount}
}}{\end{list}}

\newcounter{textItem}

\newcounter{condItem}

\ifx\undefined\slot
\newcommand{\slot}{\,-\,}
\fi
\newcommand{\lmod}[1]{{#1}_{\text{mod}}^{-1}}
\newcommand{\lalg}[1]{{#1}_{\text{alg}}^{-1}}

\newcommand{\modmap}[2]{[#1,#2]_{\text{mod}}}

\DeclareMathOperator{\sym}{Sym}
\newcommand{\modmod}{/\!/}
\newcommand{\rmod}[1]{#1\text{-mod}}
\newcommand{\ralg}[1]{#1\text{-alg}}
\DeclareMathOperator{\ho}{ho}

\newcommand{\algcond}[1]{(\romannumeral #1)$_{\text{alg}}$}
\newcommand{\modcond}[1]{(\romannumeral #1)$_{\text{mod}}$}

\begin{document}

\maketitle


\section*{Introduction}

This paper concerns a subtle issue that arises when attempting to
invert element in the homotopy groups of equivariant commutative
($E_{\infty}$) rings.  It arose for the authors
in~\cite{non_exist_of_elemen_of_kervaire}, as part of the construction
of the spectrum there denoted $\Omega$, and was called to our
attention when Justin Noel asked us, ``how do you know you can do
that?''

To set the stage, let $G$ be a finite group, and $R$ a $G$-equivariant
commutative ring, by which we mean and equivariant $E_{\infty}$ ring
spectrum.   We suppose given a subset $S\subset \pi^{G}_{\ast}R$.   

\begin{ques}
\label{ques:1}
Can one form an equivariant commutative ring $S^{-1}R$?
\end{ques}

\begin{eg}
\label{eg:1}
Take $G=\Z/2$, $R=S^{0}$, and $S=\{2 \}\subset\pi_{0}S^{0}$.   In this
case Noel's question becomes, ``is it possible to invert $2$ in
equivariant stable homotopy theory?''
\end{eg}

We will see that the questions is a little misleading as stated.  One
can in fact always invert $2$ (in fact one can invert anything), but
in order to know that the outcome is the expected one, there is a
non-trivial condition to check.

The authors would like to thank Justin Noel for bringing this issue to
our attention.  It is underlain by a subtlety in the symmetric
monoidal structure on equivariant spectra that seems to have gone
unnoticed.  It is also intimately connected to the theory of the
``norm'' in~\cite{non_exist_of_elemen_of_kervaire}.  This note
describes the main ideas that go into resolving Question~\ref{ques:1}.
Details and a much more expanded discussion will appear elsewhere.

A purely algebraic analogue of some of the structures discussed in
this paper were first formulated by Tambara in~\cite{MR1209937}, and
are now known by the name of ``Tambara functors.''  For a recent
general introduction to the algebraic theory the reader is referred
to~\cite{stricland:_tambar}.  The question of localization in the
context of Tambara functors has been investigated by
Nakaoka~\cite{nakaoka:_mackey_tambar}, in which his Proposition~4.8 is
similar to our Proposition~\ref{thm:4}.

During the preparation of this work the authors were supported by
supported by DARPA through the AFOSR grant number
HR0011-10-1-0054-DOD35CAP.  The first author was also supported by NSF
DMS--0906285, NSF DMS--1207774, and the Sloan Foundation.  The second
author received additional support from the NSF through the grant
DMS-0906194

\section{Localization of commutative rings}
\label{sec:local-comm-rings}
\numberwithin{equation}{section}

Let's begin with localization in ordinary commutative algebra.
Suppose that $R$ is a commutative ring, and $S\subset R$ a subset of
elements we intend to invert.

\begin{defin}
\label{def:1} The map $R\to S^{-1}R$ is the universal map $f:R\to R'$
having the property that $f(s)$ is a unit in $R'$ for every $s\in S$.
\end{defin}

Spelled out more precisely, the condition is that given a map $f:R\to
R'$ with the property that $f(s)$ is a unit for every $s\in S$, there
is a unique map $S^{-1}R\to R'$ factoring $f$.  

Such a universal arrow clearly exists.  One can take $S^{-1}R$ to be
the quotient
\[
S^{-1}R=R[x_{s} \mid s\in S]/(s\cdot x_{s}=1)
\]
of the polynomial algebra over $R$ constructed by adding a variable
$x_{s}$ for each $s\in S$, by the relation indicated.  There are
advantages and disadvantages to this characterization.  On one hand,
it is clear that $S^{-1}R$ exists, and that the construction can be
generalized.  On the other hand, many of the fundamental properties of
$S^{-1}R$ are not apparent.  For example, it is not clear that
$S^{-1}R$ is flat, or even that it is non-zero.

There is a second common description of $S^{-1}R$.  Suppose for
simplicity that $S$ is finite and let $\{s_{1},s_{2},\dots \}$ be a
sequence of elements of $S$ having the property that each $s\in S$
occurs infinitely often as an $s_{i}$.  The second description of
$S^{-1}R$ is that it is the colimit of the sequence
\[
\cdots\xrightarrow{s_{i}}{}
R\xrightarrow{s_{i+1}}{}R\xrightarrow{s_{i+2}}{}\cdots .
\]
From this description it is obvious that $S^{-1}R$ is flat, and that
it is non-zero as long as there is an element of $R$ which is not
annihilated by any product of powers of elements of $S$.  It is not,
however, immediately obvious from this description that $S^{-1}R$ is a
commutative ring, or even that it is independent of the choice of
sequence $\{s_{i} \}$.  These are not difficult things to verify, but
the proofs involve explicit use of the elements of $R$, and so don't
generalize well.

In fact, what this second construction actually produces is a
universal $R$-{\em module} map from $R$ to an $R$-module $M$ having the
property that for each $s\in S$, the map ``multiplication by $s$''
\[
s:M\to M
\]
is an isomorphism.   Since an $R$-algebra is an $R$-module there is a
canonical map between these two constructions, and our basic question
becomes 
\begin{ques}
\label{ques:2}
When is the canonical map from the localization of $R$ as a module to
the localization of $R$ as an algebra and isomorphism?
\end{ques}

In raising the issue we have considered only $R$, first as an
$R$-algebra and then as an $R$-module.  Our formulation will better
lend itself to generalization if we think in terms of an arbitrary
commutative $R$-algebra $A$ and an arbitrary commutative $R$-module
$M$.  For a subset $S\subset R$ we let $A\to \lalg{S}A$ be the
universal map from $A$ to a commutative $R$-algebra in which each
element of $S$ becomes a unit.  Similarly, we let $M\to \lmod{S}M$ be
the universal $R$-module map to an $R$-module on which multiplication
by each element of $S$ is an isomorphism.  The universal properties
give an evident map
\[
S_{\text{mod}}^{-1}A \to S_{\text{alg}}^{-1}A
\]
and we are interested in the conditions guaranteeing it to be an
isomorphism.

\section{Homotopy theoretic commutative algebra}
\label{sec:homot-theor-comm}

We now pursue the previous considerations in homotopy theory.  Fix a
finite group $G$, and let $R$ be a $G$-equivariant $E_{\infty}$ ring
spectrum.  We will just refer to this structure by saying that $R$ is
an ``equivariant commutative ring.''  We will also choose a finite
subset $S$ of the equivariant homotopy groups of $R$, and write $|s|$
for the ``degree'' of an element $s\in S$.  We are being deliberately
noncommittal about the term ``equivariant homotopy groups.''  When we
are led to the eventual formulation of our criteria in
\S\ref{sec:equiv-mult-clos} we will want to allow $s$ to be an
$G$-equivariant map $S^{V}\to R$, where $V$ is a virtual real
representation of $G$, and $S^{V}$ its ``one-point compactification.''
In that case the symbol $|s|$ will refer to $V$.  Nothing will be lost
at this point if the reader imagines the set $S$ to consist of
$G$-equivariant maps
\[
s:S^{n}\to R
\]
with $n\in \Z$, in which case $|s|$ is the integer $n$.   

It is a little easier to reverse the order of the previous section and
begin with the localization of modules.

\subsection{Modules} \label{sec:modules} Let $\rmod{R}$ be the
category of equivariant left $R$-modules.  We do homotopy theory in
$\rmod{R}$ by defining a map to be a weak equivalence if the
underlying map of equivariant spectra is, and employ the shorthand
notation $\modmap{M}{N}$ to denote the abelian group of maps from $M$
to $N$ in the homotopy category $\ho\rmod{R}$.

For an equivariant $R$-module $M$ and $s\in S$ let
\begin{equation}
\label{eq:3}
s:S^{|s|}\wedge M\to M
\end{equation}
be the map representing ``multiplication by $s$.''  It is given by the
composite
\[
S^{|s|}\wedge M\xrightarrow{s\wedge 1} R\wedge M \to M.
\]
Since any suspension of a module is a module, this means we will be
using the same symbol $s$ to denote any suspension of the
map~\eqref{eq:3}.  Similarly, if $M$ is an equivariant {\em right}
$R$-module we will let $\cdot s$ denote the map ``right multiplication
by $s$, given by
\[
M \wedge S^{|s|} \xrightarrow{1\wedge s} M\wedge R \to M.
\]

When $M$ is an equivariant $R$ bi-module the map $s$ is a map of right
$R$-modules and $\cdot s$ is a map of equivariant left modules.  Of
course our assumption that $R$ is commutative means that every left
module can be regarded as a bimodule, and the distinction is not
actually significant.  However since our main concern is a subtlety in
the theory of equivariant commutative rings, it seems best to use the
commutativity of $R$ sparingly.

\begin{defin}
\label{def:3} A left $R$-module $M$ is $S$-local if for each $s\in S$
the map $s:S^{|s|}\wedge M\to M$ is a weak equivalence.  Similarly a
right module is $S$-local if each of the maps $\cdot s$ is a weak
equivalence.
\end{defin}

\subsubsection{The localization conditions} \label{sec:local-cond} The
module localization of $M$ is meant to be a universal arrow from $M$
to an $S$-local $R$-module $N$.  A little homotopical algebra is
needed to state the universal property accurately.  Let
$\rmod{R}_{M/}$ be the category of equivariant $R$-modules equipped
with a map from $M$.   We define a map
\[
M\to N_{1}\to N_{2}
\]
in $\rmod{R}_{/M}$ to be a {\em weak equivalence} if the underlying
map $N_{1}\to N_{2}$ of equivariant spectra is.  We say that an object
$M\to N$ is {\em $S$-local} if $N$ is.  The module localization $M\to
\lmod{S}M$ is characterized uniquely up to unique isomorphism in
$\ho\rmod{R}_{M/}$ by the following conditions:
\begin{itemize}
\item[\modcond{1}] The module $\lmod{S}M$ is $S$-local;
\item[\modcond{2}] If $M\to N$ is
$S$-local, then $\ho\rmod{R}_{M/}(S^{-1}M,N) = \text{pt}$.
\end{itemize}
\noindent In stating these conditions we have abbreviated $M\to
\lmod{S}M$ and $M\to N$ to just $\lmod{S}M$ and $N$ in order to avoid
a cumbersome expression for maps in the homotopy category.

\subsubsection{Local modules}
\label{sec:local-modules}
To construct the $S$-localization of $M$ we first investigate
$S$-local modules.  let $\{s_{1},s_{2},\dots \}$ be a sequence of
elements of $S$ having the property that each $s\in S$ occurs
infinitely often as an $s_{i}$, and let $E$ be the right $R$-module
given by 
\[
E = \ho\varinjlim \left\{R \xrightarrow{s_{1}} S^{-|s_{1}|}\wedge R\xrightarrow{s_{2}}{}
S^{-|s_{1}|-|s_{2}|}\wedge R\to\cdots\right\}.
\]
and $E'$ left $R$-module
\[
E' = \ho\varinjlim \left\{R \xrightarrow{\cdot s_{1}} S^{-|s_{1}|}\wedge
R\xrightarrow{\cdot s_{2}}{}
S^{-|s_{1}|-|s_{2}|}\wedge R\to\cdots\right\}.
\]
Of course since $R$ is commutative there is a canonical equivalence
between $E$ regarded as an equivariant left $R$-module, and $E'$.
Write
\[
\mathcal Z =\{S^{t}\wedge (R \cup_{\cdot s}C(S^{|s|}\wedge R))\mid
s\in S, t\in\Z \}
\]
for the set of suspensions mapping cones of the maps $s$.  The
elements of $\mathcal Z$ are left $R$-modules.

Though we eventually hope to endow $E$ with the structure of an
equivariant commutative $R$-algebra, it isn't clear at this point
whether or not it has a multiplication.  It is, however, a right
$R$-module, and it does come equipped with a ``unit'' $R\to E$.

\begin{lem}
\label{thm:7} Both $E$ and $E'$ are $S$-local.  For each $Z\in\mathcal
Z$ the $G$-spectrum $E\underset{R}{\wedge}Z$ is contractible.
\end{lem}

\begin{pf}
This result makes use of the commutativity of $R$.  We begin with the
assertion that $E$ is $S$-local.  Let $s\in S$ and consider the
following commutative ladder
\[
\xymatrix{
R\wedge S^{|s|}  \ar[r]^-{s_{1}}\ar[d]^{\cdot s} 
 &S^{-|s_{1}|}\wedge R\wedge S^{|s|}   \ar[r]^-{s_{s}}\ar[d]^{\cdot s}   
 & S^{-|s_{1}|-|s_{2}|}\wedge R\wedge S^{|s|}  \ar[d]^{\cdot s} \ar[r]
 & \cdots \\
R \ar[r]_-{s_{1}}
 & S^{-|s_{1}|}\wedge R \ar[r]_-{s_{2}}        
 & S^{-|s_{1}|-|s_{2}|}\wedge R \ar[r] &\cdots 
}
\]
The induced map of homotopy colimits of the rows is right
multiplication by $s$ 
\[
\cdot s:E\wedge S^{|s|}\to E.
\]
One easily checks it to be an equivalence by considering the induced
map of equivariant homotopy groups, using the homotopy commutativity
of $R$ and the condition that each $s\in S$ occurs as an $s_{i}$
infinitely often.  The proof that $E'$ is $S$-local is similar.  The
last assertion follows from the cofibration sequence
\[
S^{t}\wedge E\wedge S^{|s|}
\xrightarrow{\cdot s}{}
S^{t}\wedge E
\to
E\underset{R}{\wedge}Z.
\]
\end{pf}
 
\begin{prop}
\label{thm:1}
For an $R$-module $N$ the following are equivalent
\begin{thmList}
\item\label{thmItem:1} For each $s\in S$ the map $s:S^{|s|}\wedge N\to N$ is a weak
equivalence;
\item\label{thmItem:2} The map $N\to E\underset{R}\wedge N$, given by smashing with the
unit, is a weak equivalence;
\item\label{thmItem:3} For each $Z$ with $E\underset{R}{\wedge}Z\sim\ast$, the group
$\modmap{Z}{N}$ is trivial;
\item\label{thmItem:4} For each $Z\in\mathcal Z$ the group $\modmap{Z}{N}$ is trivial;
\end{thmList}
\end{prop}

\begin{pf}
The assertion~\thmItemref{thmItem:1} $\implies$~\thmItemref{thmItem:2}
follows from the fact that smashing commutes with colimits,
and~\thmItemref{thmItem:2} $\implies$~\thmItemref{thmItem:3}
is a consequence of the diagram
\[
\xymatrix{
Z  \ar[r]^{f}\ar[d]  &  N  \ar[d]^{\sim} \\
E\underset{R}{\wedge}Z  \ar@<.8ex>[r]_{E\wedge f}        & E\underset{R}{\wedge}N
}
\]
which is constructed by smashing an equivariant $R$-module map $Z\to
N$ with the unit map $R\to E$.  The implication \thmItemref{thmItem:3}
$\implies$~\thmItemref{thmItem:4} is immediate from the second
assertion in Lemma~\ref{thm:7}.  Finally,~\thmItemref{thmItem:4}
$\implies$~\thmItemref{thmItem:1} follows from the fact that the map
\[
\modmap{S^{t}\wedge R}{N}
\to 
\modmap{S^{t}\wedge S^{|s|}\wedge R}{N},
\]
induced by $s:S^{|s|}\wedge R\to R$, is isomorphic to the map
\[
\pi_{t}N \to \pi_{t}S^{-|s|}\wedge N
\]
induced by $s$.
\end{pf}

From Lemma~\ref{thm:7} and the equivalence of the first two parts
statements in Proposition~\ref{thm:1} we have

\begin{cor}
\label{thm:6}
The map $R\to E'$ extends to a weak equivalence 
\[
E \to E\underset{R}{\wedge} E'.
\]
\qed
\end{cor}

As mentioned earlier, the commutativity of $R$ allows one to regard
left modules as right modules and for us to identify $E$ with $E'$.
It also makes the category of equivariant left modules into a
symmetric monoidal category.  With this in mind the equivalence in
Corollary~\ref{thm:6} can be written as an equivalence
\[
E\to E\underset{R}{\wedge}E,
\]
making it appear that the module $E$ {\em must} be an equivariant
commutative $R$-algebra, and that it must be so in a unique way.  As
we shall see, something like this is true.  But there is more to the
story.

\subsubsection{Construction of the module localization}
\label{sec:constr-module-local-1}

It follows from part~\thmItemref{thmItem:3} of
Proposition~\ref{thm:1} that an $R$-module is $S$-local if and only if
it is $E$-local in the sense of Bousfield~\cite{Bous:LocSpectra}, and
that the universal arrow
\begin{equation}
\label{eq:4}
M\to \lmod{S}M
\end{equation}
we are looking for is the Bousfield localization.  This leads to a
construction.  The Bousfield localization of any $R$-module $M$ can be
constructed as the homotopy colimit of the sequence
\[
M_{0}\to\cdots\to M_{i}\to M_{i+1}\to\cdots
\]
constructed inductively, starting with $M_{0}=M$ and, building
$M_{i+1}$ out of $M_{i}$ as a pushout
\begin{equation}
\label{eq:5}
\xymatrix{
Z_{i}  \ar[r]\ar[d]  &  CZ_{i} \ar[d] \\
M_{i}  \ar[r]       & M_{i+1}
}
\end{equation}
in which $Z_{i}$ is the wedge
\[
Z_{i}=\bigvee_{\substack{Z\to M_{i} \\ Z\in\mathcal Z}} Z
\]
and the map $Z_{i}\to M_{i}$ is the tautological map, whose restriction to
each summand is the map indicated by its index.

Let's temporarily write $M\to L_{E}M$ for the map from $M$ to the
homotopy colimit of the sequence above.  The fact that the modules
$Z\in\mathcal Z$ are compact objects implies that $L_{E}M$ satisfies
condition~\thmItemref{thmItem:4} of Proposition~\ref{thm:1} and hence
is $S$-local.  The pushout squares show that for any $M\to N$ in which
$N$ is $S$-local
\[
\ho\rmod{R}_{M/}(M_{i},N) = \text{pt}
\]
hence 
\[
\ho\rmod{R}_{M/}(L_{E}M,N) = \text{pt}.
\]
Thus $M\to L_{E}M$ satisfies conditions \modcond{1} and \modcond{2} of
\S\ref{sec:local-cond}, so can be taken to be the universal arrow $M\to
\lmod{S}M$ we were seeking.

Bousfield's construction leads to a further property of localization
that is special to {\em modules}.  Applying
$E\underset{R}{\wedge}(\slot)$ to the pushout square~\eqref{eq:5} one
finds that the map
\[
E\underset{R}{\wedge}M_{i}\to
E\underset{R}{\wedge}M_{i+1}
\]
is a weak equivalence.   It follows that the map 
\[
E\underset{R}{\wedge}M\to E\underset{R}{\wedge} L_{E}M
\]
is a weak equivalence.  This leads to another characterization of the
module localization $M\to \lmod{S}M$ (which, in Bousfield's
formulation~\cite{Bous:LocSpectra} is taken as the definition).
\begin{prop}
\label{thm:8} Let $M\to N$ be a map of equivariant left $R$-modules in which
$N$ is $S$-local.   If
\[
E\underset{R}{\wedge}M\to E\underset{R}{\wedge}N
\]
is a weak equivalence then the unique map $S^{-1}M\to N$ making 
\begin{equation}
\label{eq:6}
\xymatrix{
M  \ar[r]\ar[dr]  & S^{-1}M  \ar[d] \\
         & N
}
\end{equation}
commute is a weak equivalence.
\end{prop}

\begin{pf}
Smashing~\eqref{eq:6} with $E$ gives a diagram
\[
\xymatrix{
E\underset{R}{\wedge} M  \ar[r]\ar[dr]  & E\underset{R}{\wedge} S^{-1}M  \ar[d] \\
         & E\underset{R}{\wedge} N
}
\]
in which the top arrow is a weak equivalence by the above discussion
and the diagonal arrow is a weak equivalence by assumption.  This
means that the vertical is also a weak equivalence.  Now consider the
diagram
\[
\xymatrix{
S^{-1}M  \ar[r]\ar[d]  &  N \ar[d] \\
E\underset{R}{\wedge}S^{-1}M  \ar@<.8ex>[r]        & E\underset{R}{\wedge}N.
}
\]
We have just shown that the bottom horizontal arrow to be a weak equivalence.  The
vertical arrows are weak equivalences by part~\thmItemref{thmItem:2}
of Proposition~\ref{thm:1}.  It follows that the top arrow is a weak
equivalence.
\end{pf}

\begin{cor}
\label{thm:2} The unit map $R\to E'$ is the $S$-localization of $R$,
and the $R$-modules $E$ and $E'$ are independent of the choice of sequence
$\{s_{1},s_{2},\dots \}$.
\end{cor}

\begin{pf}
By Lemma~\ref{thm:7} and Corollary~\ref{thm:6} the map $R\to E'$
satisfies the conditions of Proposition~\ref{thm:8}. 
\end{pf}

We have now seen that the module localization of $R$, characterized by
the properties described in \S\ref{sec:local-cond}, has the homotopy
type we expect both the algebra and module the localizations to have,
namely that of
\[
E' = \ho\varinjlim \left\{R \xrightarrow{\cdot s_{1}} S^{-|s_{1}|}\wedge
R\xrightarrow{\cdot s_{2}}{}
S^{-|s_{1}|-|s_{2}|}\wedge R\to\cdots\right\}.
\]
Our next aim is to investigate what actually happens with the algebra
localization.

\subsection{Algebras}
\label{sec:algebras}

We now turn to the localization in algebras.  Let $\ralg{R}$ be the
category of $G$-equivariant $E_{\infty}$ algebras over $R$.  We do
homotopy theory in $\ralg{R}$ by defining a map to be a {\em weak
equivalence} if the underlying map of equivariant spectra is.  We will
say that an $R$-algebra is {\em $S$-local} if it is so when regarded
as a left $R$-module.  As with modules, in order to state the
universal property of algebra localization we will need to consider
the category of equivariant commutative $R$-algebras $B$ equipped with
an $R$-algebra map $A\to B$.  But such a $B$ is just an $A$-algebra,
so we will denote this category $\ralg{A}$, rather than using the
analogue $\ralg{R}_{A/}$ of our notation for modules.  We do homotopy
theory in $\ralg{A}$ by defining a map $A\to B_{1}\to B_{2}$ to be a
weak equivalence if the underlying map of $G$-spectra $B_{1}\to B_{2}$
is.

Given an $R$-algebra $A$ we wish to construct an $R$-algebra map $A\to
\lalg{S}A$ which is universal for maps from $A$ to $R$-algebras $B$
which are $S$-local.  As in the discussion of modules, the universal
property consists of the conditions 
\begin{itemize}
\item[\algcond{1}] The $R$-algebra $\lalg{S}A$ is $S$-local;
\item[\algcond{2}] If $A\to B$ is $S$-local, then the set
$\ho\ralg{A}(\lalg{S}A,B)$ consists of just one point.
\end{itemize}

Imitating the discussion in \S\ref{sec:local-comm-rings}, one can
certainly construct $\lalg{S}A$ as the pushout of
\[
\xymatrix{
\sym_{A}[r_{s}\mid s\in S]  \ar[r]\ar[d]  & A  \ar[d] \\
\sym_{A}[x_{s}\mid s\in S]  \ar[r]        & \lalg{S}A
}
\]
in which $|x_{s}|=-|s|$, $|r_{s}|=0$, the left vertical arrow sends
$r_{s}$ to $s\cdot x_{s}$ and the top arrow sends $r_{s}$ to $1$.
This construction produces an algebra which is clearly $S$-local,
however it requires a small computation in topological Andre-Quillen
cohomology to show that it satisfies property \algcond{2}.

A better approach is to imitate the construction of the Bousfield
localization described in~\S\ref{sec:modules} and construct
$\lalg{S}A$ as the homotopy colimit of the sequence
\[
\cdots\to A_{i}\to A_{i+1}\to\cdots
\]
constructed inductively, starting with $A_{0}=A$, and constructing
$A_{i+1}$ from $A_{i}$ as the (algebra) pushout of
\begin{equation}
\label{eq:7}
\xymatrix{
\sym_{R}(Z_{i})  \ar[r]\ar[d]  &  \sym_{R}(CZ_{i}) \ar[d] \\
A_{i}  \ar[r]       & A_{i+1},
}
\end{equation}
in which $Z_{i}$ is the wedge of $R$-modules
\[
Z_{i}=\bigvee_{\substack{Z\to A_{i} \\ Z\in\mathcal Z}} Z
\]
and the map $\sym(Z_{i})\to A_{i}$ is the algebra extension of the
tautological map whose restriction to each summand is indicated by its
index.  From this description the universal property is more obvious.
That $\lalg{S}A$ is $S$-local follows, as before, from the compactness
of the $Z\in\mathcal{Z}$.  The Meyer-Vietoris sequence associated to
each pushout square shows that if $B$ is an $S$-local $A$-algebra then
for each $i$ the restriction map
\[
\ho\ralg{A}(A_{i+1},B) \to 
\ho\ralg{A}(A_{i},B)
\]
is a bijection.   Since $A_{0}=A$ each of these sets consists of just
one point, and hence so does 
\[
\ho\ralg{A}(\lalg{S}A,B).
\]
Thus Bousfield's construction produces an $A$-algebra which is easily
verified to have the two properties characterizing the localization of
$A$ at $S$.   

Since $\lalg{S}A$ is $S$-local as an $R$-module, there is a unique a
module map
\begin{equation}
\label{eq:1}
\lmod{S}A\to \lalg{S}A,
\end{equation}
and our main question becomes

\begin{ques}
\label{ques:3}
Under what conditions on $S$ is~\eqref{eq:1} a weak equivalence?
\end{ques}

\noindent By Proposition~\ref{thm:8}, this is equivalent to

\begin{ques}
\label{ques:4}
Under what conditions on $S$ is the map 
\[
E\underset{R}{\wedge} A \to 
E\underset{R}{\wedge} \lalg{S}A 
\]
a weak equivalence?
\end{ques}

\noindent An analysis of the pushout squares~\eqref{eq:7} shows that this, in
turn, is equivalent to

\begin{ques}
\label{ques:5}
Under what conditions on $S$ is the map 
\[
E\underset{R}{\wedge}\sym_{R}(Z) \to E\underset{R}{\wedge}R
\]
induced by $Z\to\ast$ a weak equivalence for all $Z\in\mathcal Z$?
\end{ques}

\noindent Finally, this is equivalent to 

\begin{ques}
\label{ques:6}
Under what conditions on $S$ is 
\[
E\underset{R}{\wedge}\sym^{n}Z
\]
contractible, for all $n\ge 1$ and all $Z\in\mathcal Z$? 
\end{ques}

We will address this question in the next section.

\section{Extended powers and indexed monoidal products}
\label{sec:extend-powers-index}

\subsection{The homotopy type of symmetric powers}
\label{sec:homot-type-symm}

We have a spectrum $Z$ which is $E$-acyclic in the sense that
$E\underset{R}{\wedge}Z$ is contractible, and we wish to know whether
or not $\sym^{n}Z$ is $E$-acyclic.  To define $\sym^{n}Z$ requires
that $R$ be commutative so that the category $\rmod{R}$ acquires a
symmetric monoidal structure.  The definitions are arranged (see for
example~\cite{elmendorf97:_rings, smith00:_symmet, MR1922205})
so that $\sym^{n}Z$ has the homotopy type of
\[
(E_{G}\Sigma_{n})_{+}\underset{\Sigma^{n}}{\wedge}Z^{(n)},
\]
in which
\[
Z^{(n)} = Z\underset{R}{\wedge}\dots\underset{R}{\wedge}Z,
\]
is the iterated smash product over $R$, and $E_{G}\Sigma_{n}$ is the
$G$-equivariant analogue of the free contractible $\Sigma_{n}$-space
$E\Sigma_{n}$.  We first consider the case of the trivial group.


When $G$ is trivial one has
\[
\sym^{n}Z\sim
(E\Sigma_{n})_{+}\underset{\Sigma^{n}}{\wedge}Z^{(n)}.
\]
Since $E\Sigma_{n}$ is a homotopy colimit of free $\Sigma_{n}$-sets,
it suffices to show that if $T$ is a free $\Sigma_{n}$-set then
\[
T_{+}\underset{\Sigma_{n}}{\wedge}Z^{(n)}
\]
is $E$-acyclic.  This reduces to the case in which $T=\Sigma_{n}$, in
which case
\[
T_{+}\underset{\Sigma_{n}}{\wedge}Z^{(n)}=
(\Sigma_{n})_{+}\underset{\Sigma_{n}}{\wedge}Z^{(n)}=
Z^{(n)}.
\]
But if $Z$ is $E$-acyclic, then so is $Z^{(n)}$ by the associativity
of the smash product.    We therefore have

\begin{prop}
\label{thm:3} When $G$ is the trivial group, the map
$\lmod{S}A\to\lalg{S}A$ is a weak equivalence for any equivariant
commutative $R$-algebra $A$. \qed
\end{prop}


How does this change if $G$ is non-trivial?  The main difference is
that the $\Sigma_{n}$-space $E\Sigma_{n}$ must be replaced by the
$E_{G}\Sigma_{n}$.  As mentioned above, $E_{G}\Sigma_{n}$ is the
$G$-analogue of $E\Sigma_{n}$.  It is a space with a
$G\times\Sigma_{n}$-action and can be defined to be the total space of
the universal $G$-equivariant principal $\Sigma_{n}$-bundle.  It is
characterized up to $G\times \Sigma_{n}$-equivariant weak equivalence
by the property that its fixed point space for $H\subset
G\times\Sigma_{n}$ is empty if $H\cap \Sigma_{n}$ is non-zero, and
contractible otherwise.  This characterization shows that the
$E_{G}\Sigma_{n}$ is a homotopy colimit of transitive
$G\times\Sigma_{n}$-sets $T$ which are $\Sigma_{n}$-free.  Thus to
show that
\[
\sym^{n}Z = (E_{G}\Sigma_{n})_{+}\underset{\Sigma^{n}}{\wedge}Z^{(n)}
\]
is $E$-acyclic when $Z$ is, it suffices to show that for such $T$, if
$Z$ is $E$-acyclic, so is the result of forming
\[
T_{+}\underset{\Sigma_{n}}{\wedge}Z^{(n)}.
\]
We have therefore reduced Question~\ref{ques:3} to
\begin{ques}
\label{ques:7}
What conditions on $S$ guarantee that 
\[
T_{+}\underset{\Sigma_{n}}{\wedge}Z^{(n)}
\]
is $E$-acyclic when $Z$ is $E$-acyclic, where $T$ is a
$\Sigma_{n}$-free finite $G\times\Sigma_{n}$-set?
\end{ques}

\subsection{Indexed monoidal products} 
\label{sec:index-mono-prod}

To get a feeling for the construction
\[
T_{+}\underset{\Sigma_{n}}{\wedge}Z^{(n)}
\]
let's focus on a simple example.  Take $n=2$ and $G=\Z/2$ and write
$g\in G$ for the generator.   We consider two examples of $\Sigma_{2}$-free, transitive
$G\times\Sigma_{2}$-sets
\begin{align*}
T_{\Delta} &= G\times\Sigma_{2}/\text{diagonal} \quad\text{and} \\
T_{0} &= G\times\Sigma_{2}/G.
\end{align*}
For a $G$-spectrum $Z$ one easily checks that 
\[
(T_{0})_{+}\underset{\Sigma_{2}}{\wedge}Z\wedge Z
\]
is just $Z\wedge Z$ with the diagonal $G$-action 
\[
x\wedge y\mapsto g(x)\wedge g(y)
\]
and 
\[
(T_{\Delta})_{+}\underset{\Sigma_{2}}{\wedge}Z\wedge Z
\]
is $Z\wedge Z$ with $G$-action 
\[
x\wedge y \mapsto g(y)\wedge g(x).
\]
For instance, let Let $Z=G_{+}$, i.e. $S^{0}\vee S^{0}$ with the
``flip'' action of $G$.  Then the two constructions work out as
indicated below: 

\bigskip
\begin{center}
\includegraphics[]{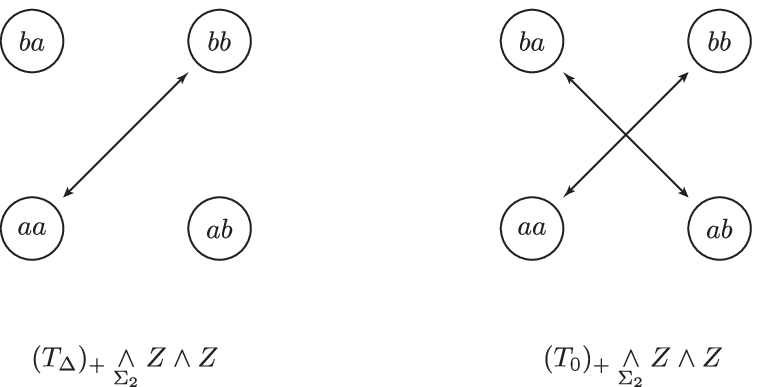}
\end{center}
\bigskip

\noindent and when $Z$ is the suspension spectrum of $\{a,b \}_{+}$ with the 
trivial action of $G$ one gets

\bigskip
\begin{center}
\includegraphics[]{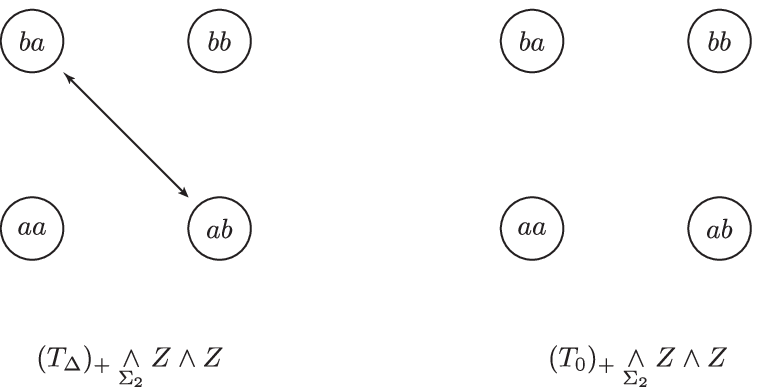}
\end{center}
\bigskip

In the general case $T_{+}\underset{\Sigma_{n}}{\wedge}Z^{(n)}$
is constructed as a combination of the following processes.  First,
starting with a subgroup $H\subseteq G$, an $H$-spectrum $X$ and a finite
$H$-set $U$ one forms
\[
W= \bigwedge_{u\in U}X.
\]
This is a spectrum with an $H$-action.   Next one induces the result
up to a $G$-spectrum by forming
\[
G_{+}\underset{H}{\wedge}W = \bigvee_{G/H} g\cdot W.
\]
These processes are studied in much more detail
in~\cite{non_exist_of_elemen_of_kervaire} where they are called
``indexed monoidal products.''  Their similarity is more apparent with
a slightly more general setup.  Let $U$ be finite $H$-set and
$u\mapsto X_{u}$ a functor from the action category $U\modmod H$ to
spectra.  Thus $X_{u}$ is equivariant for the isotropy group
$H_{u}\subset G$ of $u$, and there are maps $g_{\ast}:X_{u}\to X_{gu}$
whose formation is associative, etc.  Out of this functor one can
construct $H$-spectra by forming the ``indexed sum''
\[
\bigvee_{u\in U} X_{u},
\]
the ``indexed product''
\[
\prod_{u\in U} X_{u},
\]
and the ``indexed smash product''
\[
\bigwedge_{u\in U} X_{u}.
\]

Pursuing this, what emerges is the principle that in equivariant homotopy theory it
isn't enough to just study ordinary coproducts, products and smash products.
One must also systematically include the situation in which the group
$G$ is acting on the indexing set.  The benefits from doing this
impact even the most basic notions of the field.  For example, the
additivity of non-equivariant stable homotopy theory is expressed by
the assertion that for any $X, Y$, the map
\[
X\vee Y\to X\times Y
\]
is a weak equivalence.   Equivariant stable homotopy theory possesses
an  enhanced {\em equivariant additivity}:  the map 
\[
\bigvee_{u\in U} X_{u}\to
\prod_{u\in U} X_{u}
\]
is a weak equivalence.  In the presence of stability (the fact that
fibration sequences are weakly equivalent to cofibration sequences),
this property is equivalent to the dualizability of finite $G$
CW-spectra and to the invertibility of representation spheres $S^{V}$.
The invertibility of representation spheres is arguably {\em the} 
fundamental characteristic of equivariant stable homotopy theory.  It
is often referred to as {\em equivariant stability}, and is the only
aspect of equivariant stable homotopy theory that is not easily
described in algebraic terms.  But (for finite $G$) once one thinks in
terms of indexed coproducts, equivariant stability becomes simply the
enhancement of ordinary stability one finds by requiring {\em equivariant
additivity}.  And that is a fundamentally algebraic notion.

In the next section we will turn to the structures that arise from
systematically considering equivariant multiplicativity.

\section{Equivariant multiplicative closure}
\label{sec:equiv-mult-clos}

The indexed monoidal smash product allows one to express
an enhanced commutativity property possessed by equivariant
commutative rings.  For a non-equivariant commutative ring, the
multiplication map 
\[
R\wedge R\to R
\]
is a ring homomorphism.  This is, of course, related to the fact that
the coproduct of $E_{\infty}$ rings is weakly equivalent to the smash
product.  One might naturally guess that the indexed coproduct of
equivariant commutative rings is weakly equivalent to the smash
product
\[
\coprod_{u\in U} R_{u} \approx \bigwedge_{u\in U} R_{u},
\]
and that there is canonical equivariant ring homomorphism
\begin{equation}
\label{eq:2}
\bigwedge_{u\in U} R_{u} \to R.
\end{equation}
This is in fact the case.  The particular proof depends on the
particular model for equivariant stable homotopy theory you are using,
but ultimately it comes down to the fact that the indexed product of
the universal $G_{u}\times \Sigma_{n}$-spaces
\[
\prod_{u\in U} E_{G_{u}}\Sigma_{n}
\]
is a model for the universal $G\times\Sigma_{n}$-space
$E_{G}\Sigma_{n}$.  For more details see of~\cite[\S2,
App.~A,~B]{non_exist_of_elemen_of_kervaire}.  The same discussion
applies to an equivariant commutative $R$-algebra $A$, and with the
smash product replaced with the smash product over $R$.

These considerations lead to two related notions of {\em equivariant
multiplicative closure}.

\begin{defin}
\label{def:4} A collection $\mathcal C$ of pairs $(Z,H)$ with $H$ a
subgroup of $G$ and $Z$ and $H$-equivariant $R$-module is {\em
equivariantly multiplicatively closed} if it is stable under the
formation of indexed smash products.  It is an {\em equivariant ideal}
if for every finite $H$-set $U$ and every functor $Z$ from $U\modmod
H$ to spectra, if for some $u$, $(Z_{u},H_{u})$ is in $\mathcal C$
then so is the indexed smash product
\[
\bigwedge_{u\in U} Z_{u}.
\]
\end{defin}

\begin{eg}
\label{eg:3}
If $A$ is an equivariant commutative $R$-algebra,  the set
\[
\{(Z,H)\mid A\underset{R}{\wedge}Z\underset{H}{\sim}\ast \}
\]
 of {\em $A$-acylic $R$-modules} is an equivariant ideal.  (The symbol
$\sim_{H}$ indicates an $H$-equivariant weak equivalence.)  Indeed if
for some $u$
\[
A\underset{R}{\wedge}Z_{u}
\]
is weakly $H_{u}-$contractible, then the indexed smash product
\[
\bigwedge_{u\in U} \left(A\wedge Z_{u} \right)
\]
is $H$-contractible.  The indexed smash product is a symmetric
monoidal functor, so this implies that
\[
\left(\bigwedge_{u\in U}A \right)\wedge \left(\bigwedge_{u\in U}Z_{u}\right)
\]
is $H$-contractible.   But then so is
\[
(A)\wedge \left(\bigwedge_{u\in U}Z_{u}\right)
\]
since there is an equivariant ring map 
\[
\bigwedge_{u\in U} A\to A.
\]
\end{eg}

Example~\ref{eg:3} bears directly on the situation we found ourselves
considering in the previous section.  In fact it tells us that the map
\[
\lmod{S}A\to \lalg{S}A
\]
is a weak equivalence for every $A$ if and only if it is so for $R$,
and that this holds if the set of $E$-acyclic modules forms an ideal.
By the above discussion, this condition is also necessary, since after
all, what we're aiming to show in the first place is that $E$ is an
equivariant commutative ring.

We can address this issue entirely in terms of $S$ by considering a
second notion of equivariant multiplicative closure.  Let $V$ be a
$G$-equivariant vector bundle over $U$ and $S^{V}$ its one-point
compactification.  Writing $V_{u}$ for the fiber over $u$ we may think
of $S^{V}$ as the indexed wedge
\[
\bigvee_{u\in U} S^{V_{u}} .
\]
A $G$-equivariant map
\[
\bigvee_{u\in U} S^{V_{u}} \to R
\]
corresponds to a collection of compatible $G_{u}$-equivariant maps 
\[
\alpha_{u}:S^{V_{u}}\to R.
\]
Forming the indexed smash product and composing with~\eqref{eq:2}
gives a map
\[
\bigwedge_{u\in U}S^{V_{u}} \to 
\bigwedge_{u\in U}R  \to R
\]
which we call the {\em indexed product} of the $\alpha_{u}$.   We
denote this indexed product by
\[
\prod_{u\in U}\alpha_{u}:S^{W}\to  R
\]
in which
\[
W=\bigoplus_{u\in U} V_{u}
\]
is the space of global sections of $V$.   

\begin{eg}
\label{eg:4} Consider the case $R=S^{0}$.  Then $\pi_{0}^{G}R$ is the
Burnside ring $A(G)$ of finite $G$-sets.  A functor $T(u)=T_{u}$
from $U\modmod G$ to the category of finite sets determines a
$G$-equivariant homotopy class of maps
\[
\bigvee_{u\in U} S^{0} \to S^{0}.
\] 
The indexed product in this case corresponds to the $G$-set given by
the indexed product
\[
\prod_{u\in U} T_{u}.
\]
The examples at the end of \S\ref{sec:index-mono-prod} illustrate this
situation in the case $G=\Z/2$. 
\end{eg}

\begin{eg}
\label{eg:5} Consider a general finite $G$ and a transitive $G$-set
$U=G/H$ for some proper subgroup $H\subset G$.  The indexed product of
set of $G_{u}$-equivariant maps
\[
\alpha_{u}:S^{1}\to R
\]
will be a $G$-equivariant map 
\[
\prod_{u\in U}\alpha_{u} :S^{W}\to R
\]
where $W$ is the real $G$-representation induced from the trivial
representation of $H$.  
\end{eg}

Example~\ref{eg:5} shows that even if one were only interested in
$\pi_{n}^{G}R$ for $n\in \Z$, the presence of indexed products for
equivariant commutative rings would necessitate consideration of the
groups $\pi_{V}^{G}R$ for certain virtual representations $V$ of $G$.
In fact by the time everything is all laid out, one is actually led to
consider all of the homotopy groups $\pi^{H}_{V}R$ with $H\subset G$ a
subgroup, and $V$ a virtual representation of $H$.  In order not to
overly complicate the discussion we have restricted our focus in this
groups $\pi_{V}^{G}R$.  In view of this, we henceforth take the
expression ``equivariant stable homotopy groups'' to mean the sum of
the groups $\pi_{V}^{G}$, and view it as graded over the real
representation ring of $G$.  Elements of the equivariant stable
homotopy groups will always be assumed to be homogeneous.

The enhancement of equivariant multiplication given by the formation
of ``indexed products'' goes back to the work of Greenlees and
May~\cite{MR1491447}, and is related to concept of Tambara functors in
the work of Brun~\cite{MR2308231}.  It leads to our second notion of
{\em equivariant multiplicative closure}.

\begin{defin}
\label{def:5} Let $R$ be an equivariant commutative ring.  A subset of
the equivariant stable homotopy groups of $R$ is {\em multiplicatively
closed} if it is closed under the formation of indexed products.
\end{defin}

\begin{eg}
\label{eg:2} Let $A$ be an equivariant commutative ring.  The set of
$\alpha\in\pi_{V}^{G}R$ which are {\em units} is multiplicatively
closed.
\end{eg}

Example~\ref{eg:2} tells us that if we are going to force all the
elements of some set $S$ to become units in an equivariant commutative
ring, we are also going to force all indexed products of those
elements to become units as well.   This is the subtle difference
between the equivariant and non-equivariant case.  

\begin{defin}
\label{def:2} Let $R$ be an equivariant commutative ring, and $S$ a
subset of the equivariant homotopy groups of $R$.  The {\em
equivariant multiplicative closure} $\hat S$ of $S$ is the smallest
multiplicatively closed subset of the equivariant homotopy groups of
$R$ containing $S$.
\end{defin}

We then have

\begin{prop}
\label{thm:4} Let $R$ be an equivariant commutative ring and $A$ an
equivariant commutative $R$-algebra.  For every subset $S$ of the
equivariant homotopy groups of $R$, the maps
\[
\lalg{S}A\to 
\lalg{\hat S}A
\]
and
\[
\lalg{\hat S}A\to \lmod{\hat S}A
\]
are weak equivalences.    \qed
\end{prop}

\begin{prop}
\label{thm:5} With the above notation, the map
\[
\lmod{S}A\to \lmod{\hat S}A
\]
is a weak equivalence if and only if every indexed product
of elements of $S$ divides an ordinary product of elements of $S$.
\end{prop}

Combining Propositions~\ref{thm:4} and~\ref{thm:5}, and specializing
to the case of interest, we conclude
\begin{cor}
\label{thm:9}
The map
\[
\lmod{S}R\to \lalg{S}R
\]
is a weak equivalence if and only if every indexed product
of elements of $S$ divides an ordinary product of elements of $S$. \qed
\end{cor}

The following example will be discussed in much greater detail in a
later paper.  
\begin{eg}
\label{eg:6} Consider the case $R=S^{0}$ and $G=\Z/2$.  Let $\sigma$
be the $1$-dimensional sign representation and
$s\in\pi_{-\sigma}S^{0}$ the map obtained by desuspending the fixed
point inclusion $S^{0}\to S^{\sigma}$.  Set $S=\{s \}$.  Then
$\lmod{S}R$ is the suspension spectrum of the reduced suspension of
$E\Z/2$, and is not (equivariantly) contractible.  On the other hand
the indexed product
\[
(T_{\Delta})_{+}\underset{\Sigma_{2}}{\wedge} S^{-\sigma}\wedge
S^{-\sigma} \to R
\]
depends only on the non-equivariant map underlying $s$, which is null.
It follows that the equivariant multiplicative closure of $S$ contains
$0$ and so $\lalg{S}R\sim\ast$.  Thus in this case the map
$\lmod{S}R\to \lalg{S}R$ is not a weak equivalence.
\end{eg}

\subsection{Inverting $2$} \label{sec:inverting-2} We can now turn to
our specific example of Noel's question.  Take $G=\Z/2$ and $R=S^{0}$
and $S=\{2 \}$.  As mentioned in Example~\ref{eg:4},
$\pi_{0}^{G}S^{0}$ is the Burnside ring $A(G)$.  When $G=\Z/2$ it is
one has
\[
A(G) \approx \Z[\rho]/\rho^{2}-2\,\rho
\]
where $\rho=G$ is the free, transitive $G$-set.  Under the map
$A(G)\to\pi_{0}^{G}S^{0}$, the element $2$ corresponds to the $G$-set
with $2$ elements and trivial $G$-action.  We calculated the
non-trivial indexed product of this set with itself at the end of
\S\ref{sec:index-mono-prod}.  A glance at the illustrations gives the
formula
\[
\prod_{u\in \Z/2} 2 = 2+\rho.
\]
Noel's question is thus reduced to

\begin{ques}
\label{ques:1}
Does $2+\rho$ divide a power of $2$ in $Z[\rho]/(\rho^{2}-2\,\rho)$?
\end{ques}

If the answer to the above question is ``yes'' then inverting $2$ in
an equivariant commutative ring has the expected effect.   If it
doesn't then there is something else going on.   Fortunately the
answer to this question is ``yes'' 
\[
(2+\rho)(4-\rho)=8.
\]
In fact, for any finite group $G$ inverting any integer $n$ has the
expected effect.  There is, however, a condition to check.

\bibliographystyle{amsplain}

\def\cprime{$'$}
\providecommand{\bysame}{\leavevmode\hbox to3em{\hrulefill}\thinspace}
\providecommand{\MR}{\relax\ifhmode\unskip\space\fi MR }
\providecommand{\MRhref}[2]{%
  \href{http://www.ams.org/mathscinet-getitem?mr=#1}{#2}
}
\providecommand{\href}[2]{#2}

\end{document}